# Masonry elements strengthened through Textile-Reinforced Mortar: application of detailed level modelling with a free open-source Finite-Element code

Ingrid Boem, PhD         OCRID: 0000-0003-3646-0642         e-mail: ingrid.boem@fsv.cvut.cz

Department of Concrete and Masonry Structures, Faculty of Civil Engineering, Czech Technical University, Thákurova 7, 166 29 Praha 6, Czech Republic

**Abstract:** The paper concerns the modelling of masonry elements strengthened through Textile Reinforced Mortar (TRM), a near surface system made of fiber-based grids or textiles embedded in mortar layers. Recently, the author, focusing on the mechanical characterization of TRM composites, developed a detailed level modelling approach by using the free, open-source Finite-Element code OOFEM, for the simulation of experimental tests on TRM coupons (pull-out tests, tensile tests, shear bond tests and in-plane shear tests). The model was capable to account for the failure of single components (e.g. the fibers tensile failure, the mortar cracking and crushing), as well as of their interactions (the debonding of the fibers from the mortar and of the mortar from the masonry substrate). In this paper, the detailed-level modelling approach is applied to the simulation of TRM strengthened masonry elements subjected to diagonal compression, in-plane and out-of-plane bending tests, investigating on the typical failure modes of masonry. Non-linear static analyses are performed, with nonlinearities of materials and interfaces deduced from experimental evidences. The comparison with some experimental results and a parametric study allowed to evidence the reliability of the models and their sensitivity to the main components characteristics.

## 1. Introduction

The strengthening of existing masonry through Textile-Reinforced Mortar (TRM) consists in the application of mortar layers with embedded grids or textiles made of long fibers (e.g. glass, carbon, basalt…). Different combinations of inorganic matrix and reinforcement reflect in different TRM layouts, sometimes categorized, e.g., as FRM (Fabric/Fiber Reinforced Mortar), FRCM (Fiber-Reinforced Cementitious Matrix/Mortar), TRC (Textile Reinforced Concrete), IMG (Inorganic Matrix-Grid composites), CRM (Composite Reinforced Mortar)... The application of TRM can be an effective strengthening strategy against many recurrent damage modes affecting historic, unstrengthened masonry buildings subjected to exceptional lateral actions, such as earthquakes. The benefits of several TRM systems on the mechanical performances of masonry elements against the typical failure mechanism related to in-plane actions (diagonal cracking and in-plane bending failure mechanisms) and out-of-plane actions (out-of-plane bending failure mechanism) have already been extensively investigated experimentally, as pointed out the broad literature review recently performed by the author [1]. The influence of TRM on the masonry failure modes is typically evaluated by comparing the performances of unstrengthened and strengthened samples. The benefits of TRM basically relies on its high tensile strength-to-weight ratio, able to provide some "pseudo-ductility" to the masonry, since the reinforcement opposes to the widening of cracks [2]-[3]. The variety of setups revealed the burden of testing exhaustively the reinforcement effectiveness, since strictly related to the selected configurations in terms of materials, mutual interactions, geometry, load pattern and boundary conditions. Clearly, experimental tests alone do not allow to cover the whole variety of possible configurations and loading conditions in buildings. Numerical simulations permit to extend the experimental evidences to a wider variety of TRM layouts and setups and to deeply investigate on optimized intervention and design strategies with a reduced experimental effort, which can be limited to more simple tests for the calibration of materials and interfaces laws and to target validation tests. However, the literature analysis of the numerical studies available on TRM strengthened masonry elements, performed in [1], evidenced the lack of a comprehensive approach, rather than models calibrated and applied for the reproduction of a specific test setup and combination of materials.

In this context, the EU-funded project "conFiRMa" [4]-[5], acronym for "Calibration of a numerical model for Fiber-Reinforced Mortar analysis", is aimed at calibrating and validating a numerical method, implemented in the free, open source finite element code OOFEM [6], for the assessment of the structural performances of historical masonry buildings strengthened with TRM. The purpose is to develop a multi-level approach, starting with the detailed modelling of components (masonry, mortar matrix, fiber-based reinforcement, mutual interfaces), followed by an optimization procedure to get a computationally efficient intermediate level model (by using layered elements) for the calibration of the lumped plasticity model for the global analysis of structures.

Firstly, the effort focused on the mechanical characterization of TRM [7], with a deep analysis of the state of the art on the topic and the calibration and validation of the detailed level modelling on the basis of some experimental tests on TRM coupons available in the literature (pull-out tests, tensile tests, shear bond tests and in-plane shear tests). As a continuation of the study, the present paper deals with the application of the detailed level modelling on TRM strengthened masonry elements. In particular, the detailed level modelling approach is applied to the simulations of TRM strengthen masonry elements subjected to diagonal compression tests, in-plane and out-of-plane bending tests. The aim is to prove, through comparison with experimental results, the capability of the model in reproducing the typical failure modes of masonry elements, also at the varying of the main components characteristics and test setup. The numerical study herein presented has the purpose to adopt the same modelling assumptions and parameters (independently calibrated using TRM characterization tests and tests on plain masonry elements) to provide reliable simulations under different loading conditions (in-plane and out-of-plane), as well as at the varying of the masonry and strengthening system characteristics, explicitly considering the possible damage that can affect both individual materials and interfaces, also at advanced damage levels and when substantial variations in the components are considered.

In the paper are reported: the description of the investigated TRM technique and of the main experimental outcomes, the description the numerical model features, the analysis of the numerical results, the comparison with the experimental tests and a sensitivity analysis.

## 2. The CRM strengthening technique

The TRM reinforcement considered in this study, specifically classified as Composite Reinforced Mortar (CRM), consists in the application, on the masonry wall surfaces, of a mortar layer of a nominal thickness of 30 mm, with Alkali-Resistant Glass Fiber-Reinforced Polymer (GFRP) meshes embedded (Fig. 1a). GFRP L-shape passing-through connectors, injected into holes drilled the masonry and provided with GFRP mesh devices, are also applied to improve collaboration with the substrate. The mesh is produced by twisting the yarns in the warp direction across those in the weft one, which fibers remain parallel; a thermo-hardening epoxy resin (vinylester-based) provides the yarns coating. Typically, the mesh has a squared grid pitch of 66x66 $mm^2$ and the dry fibers area in each yarn is 3.8 $mm^2$ (type "S"). Meshes with 33x33 $mm^2$ or 99x99 $mm^2$ are also available, but are less frequently used. The main mechanical characteristics of the GFRP yarns and intersections are resumed in Table 1.

CRM systems distinguish from other TRM systems, such as FRCM, for the thicker mortar coating (30-50 mm) and the type of embedded reinforcement (fiber-based open meshes with pitch ≥30 mm, pre-formed through coating with organic resin). These peculiarities make CRM systems particularly suitable for the application on historical masonry, since the typical irregularities of the rough masonry surfaces can be easily levelled within the thickness of the mortar layer. Moreover, the installation is rather simple, since, once the pre-formed composite grid has been positioned, the mortar matrix can be easily applied/sprayed, as a common plaster, without the necessity of a two-stages application and in-situ lay up of adhesive promoter, as, e.g., for FRCM systems. CRM systems are particularly prone to the use of low-medium performance mortars, which have greater chemical (ensuring breathability) and mechanical (having comparable stiffness) compatibility with historical masonry, in respect to high-performance, cement-based mortars, typical of FRCM systems.

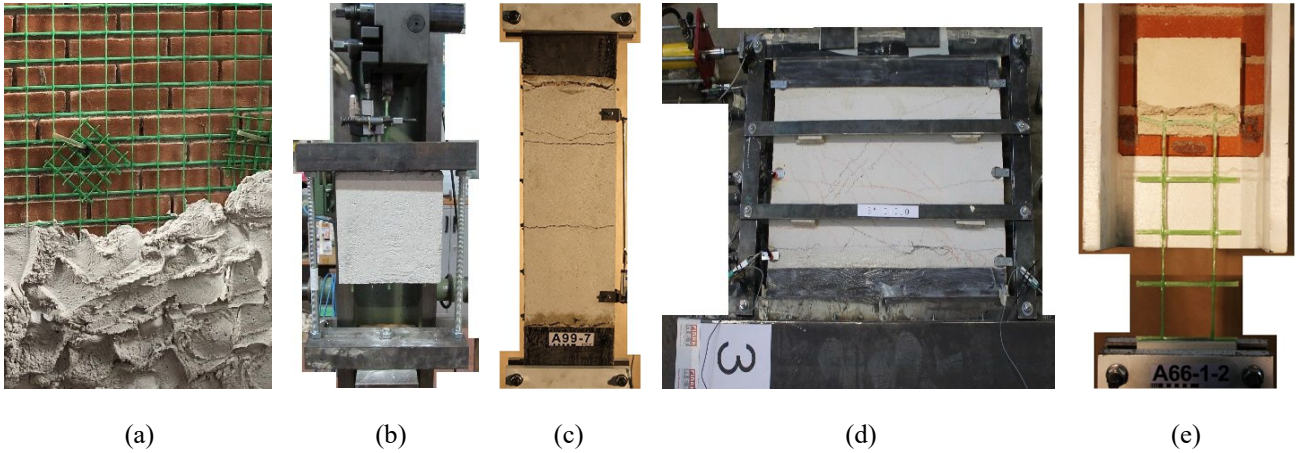

| (a) | (b) | (c) | (d) | (e) |

**Fig. 1.** The CRM strengthening technique: (a) application and some examples of characterization experimental tests such as (b) pull-out tests, (c) tensile tests, (d) in-plane shear tests and (d) shear-bond tests.

**Table 1**. Mean mechanical characteristics of GFRP yarns and intersections ($A_{w,fib}$ is the dry fibers cross section, $T_w$ the tensile resistance, $E_w$, the Young modulus, $K$ the resistance of the mesh intersection).

| Glass-fibers wires | $A_{w,fib}$ [mm²] | $T_w$ [kN] | $E_w$ [GPa] | $K$ [kN] |
|---|---|---|---|---|
| Parallel fiber yarns, type S | 3.8 | 5.62 | 69.5 | 0.55 |
| Twisted fiber yarns, type S | 3.8 | 4.49 | 62.9 | 0.48 |

## 3. Resume of experimental evidences

Many experimental results are available in the literature for the considered technique, investigating both on the CRM characterization and on its influence on the masonry performances; there is thus a wide database of experimental outcomes that can serve for the calibration and validation of a robust numerical model.

In particular, in regard to the composite characterization, the performed experiments [8]-[9] concerned tensile tests on GFRP yarns, shear tests of on mesh intersections, pull-out tests of yarns from the mortar coating (Fig. 1b), tensile tests (Fig. 1c) and in-plane shear tests (Fig. 1d) on CRM coupons and shear-bond tests on CRM layers applied on masonry wallets (Fig. 1e). The main features of these tests and the main results have already been resumed and commented in a recent paper [7], to which it is remanded for further insights. However, it's important to point out that different types of CRM failures were detected, namely the tensile rupture of yarns, the debonding of yarns from the mortar matrix and the debonding of the mortar from the masonry substrate. These occurrences depend on multiple factors, such as the reinforcement ratio, the pitch and orientation of the mesh, the bonding length, the roughness of the support...

In regard to the role of CRM in the resistant mechanism of masonry elements, experimental diagonal compression tests, in-plane bending tests and out-of-plane bending tests were performed. It is useful to summarize herein the main features and results of the tests, to allow comparison with the numerical model outcomes.

### 3.1. Diagonal compression tests

Several experimental diagonal compression tests on CRM strengthened masonry samples were performed and described in the literature [10]-[12], concerning square masonry panels (side 1160 mm) of different masonry types and thickness, such as solid brick masonry (250 and 380 mm thick), rubble stone masonry (400 and 700 mm thick) and cobblestone masonry (400 mm thick). Also different types of mortar for the coating were considered and unreinforced samples were also tested, as reference. To perform the tests, metallic stiff devices were installed at two opposite corners of the sample, to apply loading-unloading cycles along one diagonal (Fig. 2a).

The main experimental samples features and results herein considered for comparison with the experimental outcomes are summarized in Table 2: the id. starts with the indication of the type of test, DC, and masonry (e.g. B, followed with suffix "m" in case of increased thickness). Then, for the strengthened samples, there is the indication of the pitch and type of GFRP mesh and of the mortar coating (the different mortar types used for the coating have been distinguished by label "C$x$", being $x$ its approximate compressive strength). Differently, unreinforced masonry samples are distinguished with suffix U. The capacity curves are reported in dashed lines in Fig. 3, in terms of applied load, $F_{DC}$, at the varying of the shear strain, $\gamma_{DC} = \varepsilon_t + |\varepsilon_c|$ ($\varepsilon_t$ and $\varepsilon_c$ were the tensile and compressive strains monitored along the sample diagonals, measured on a base length of about 1075 mm). Generally, at least two samples were tested for each configuration: the mean values of peak resistance ($F_{DCmax}$) and respective shear strain $\gamma(F_{DCmax})$ are resumed in Table 2. Moreover, it is reported also the load in correspondance of a shear strain $\gamma_{DC} = 0.8\%$, identified as $F_{DC8}$, which provides an indication of the pseudo-ductility of the panel, intended as the capacity to preserve resistance after damage. The ratios between the peak loads of the strengthened and unstrengthened configurations ($F_{DCmax(R)}/F_{DCmax(U)}$) and between $F_{DC8}$ and $F_{DCmax}$ are also calculated.

The unreinforced solid brick masonry samples experienced a rapid drop in resistance after the sudden occurrence of the diagonal cracking, which involved mainly the mortar joints (Fig. 2b); the residual load was related to friction between elements across the cracks and to the compressive strength of the masonry struts. The unreinforced cobblestone specimens (Fig. 2c) performed a lower peak load resistance but the decrease was smoother, due to interlocking among units, providing some contrast to the crack opening. Rubble stone unreinforced masonry had an intermediate behaviour. In general, the diagonal trend of the cracks, in the center of the panel, tended to deviate at the extremities towards the ends of the steel devices. The $F_{DC8}/F_{DCmax}$ ratios resulted very low for solid brick masonry (0.16-0.34) and tended to increase for the stone samples (0.58-0.79), due to interlocking.

In strengthened samples (Fig. 2d), a diagonal crack firstly formed in the mortar coating in the center of the panel. Then, the cracking zone gradually spread, also with the formation of other cracks. Cracking extended also to the masonry, with a trend similar to that evidenced in unreinforced specimens. The GFRP mesh intervened in the cracked areas of the mortar coating, by contrasting the opening of cracks. This resulted in a gradual decrease of the resistance, as clearly emerged by analyzing the $F_{DC8}/F_{DCmax}$ ratios in strengthened samples and, in particular, by comparing DC-B-66S-C8 (0.83) and DC-B-00-C8 (0.25) - this latter had the bare plaster. The progressive collapse of several GFRP wires in the widely damaged area occurred from values of tensile deformation $\gamma_{DC}$ of about 0.8-0.9%; diffuse damage of the mortar among cracks (with local detachments of mortar portions covering the GFRP mesh) was also observed at elevate shear strains. The CRM effectiveness, in terms of resistance (ratio $F_{DCmax(R)}/F_{DCmax(U)}$), generally resulted higher for the weaker masonry. Moreover, looking at DC-B-66S samples, it emerged that a stronger mortar for the coating can lead to higher resistances, but the load decrease after the peak is faster (decrease of the $F_{DC8}/F_{DCmax}$ ratio). It should be noted that in solid brick strengthened samples, the peak load was generally associated to the formation of the first crack and subsequent load decrease, while, in the most of the stone samples, a further, even slight, load increase was monitored after this occurrence.

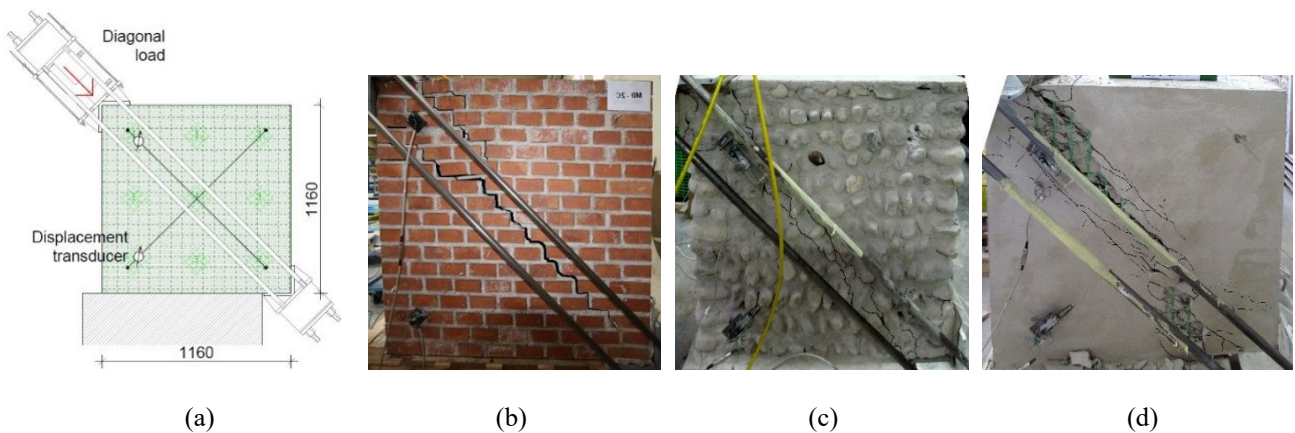

(a) (b) (c) (d)

**Fig. 2** Experimental diagonal-compression tests: (a) test setup and crack pattern unstrengthened samples (b) DC-B-U and (c) DC-C-U and of (d) strengthened sample DC-B-66S-C8.

**Table 2**. Experimental diagonal-compression tests ($F_{DCmax}$ is the peak load and $\gamma(F_{DCmax})$ the shear strain at peak load, $F_{DC8}$ the load in correspondance of a shear strain $\gamma_{DC} = 0.8\%$).

| ID | Masonry type | Masonry thickness [mm] | Mortar Type* | GFRP mesh type | $F_{DCmax}$ [kN] | $\gamma(F_{DCmax})$ [$^0/_{00}$] | $F_{DC8}$ [kN] | $F_{DCmax(R)}/F_{DCmax(U)}$ | $F_{DC8}/F_{DCmax}$ |
|---|---|---|---|---|---|---|---|---|---|
| DC-B-U | | | - | - | 217.3 | 0.26 | 54.66 | - | 0.25 |
| DC-B-00-C8 | | | C8 | - | 350.6 | 0.29 | 88.38 | 1.6 | 0.25 |
| DC-B-66S-C8 | B: Solid brick, | 250 | C8 | 66x66S | 324.2 | 0.32 | 269.27 | 1.5 | 0.83 |
| DC-B-66S-C6 | $f_{c,mortar}$~3 MPa | | C6 | 66x66S | 396.1 | 0.40 | 270.06 | 1.8 | 0.68 |
| DC-B-66S-C10 | | | C10 | 66x66S | 421.8 | 0.27 | 320.93 | 1.9 | 0.76 |
| DC-B-66S-C25 | | | C25 | 66x66S | 588.8 | 0.43 | 383.39 | 2.7 | 0.65 |
| DC-B2-U | B2: Solid brick, | 250 | - | - | 257.0 | 0.22 | 41.01 | - | 0.16 |
| DC-B2-66S-C6 | $f_{c,mortar}$~5.3 MPa | | C6 | 66x66S | 490.3 | 0.53 | 328.81 | 1.9 | 0.67 |
| DC-Bm-U | B: Solid brick, | 380 | - | - | 285.7 | 0.24 | 96.15 | - | 0.34 |
| DC-Bm-66S-C6 | $f_{c,mortar}$~3 MPa | | C6 | 66x66S | 482.5 | 0.47 | 275.66 | 1.7 | 0.57 |
| DC-R-U | R: Rubble stone, | 400 | - | - | 237.5 | 0.39 | 152.75 | - | 0.64 |
| DC-R-66S-C8 | $f_{c,mortar}$~4 MPa | | C8 | 66x66S | 461.1 | 0.89 | 404.1 | 1.9 | 0.88 |
| DC-R2-U | R2: Rubble stone, | 400 | - | - | 131.0 | 0.31 | 76.57 | - | 0.58 |
| DC-R2-66S-C6 | $f_{c,mortar}$~3 MPa | | C6 | 66x66S | 349.1 | 1.94 | 338.42 | 2.7 | 0.97 |
| DC-Rm-U | R: Rubble stone, | 700 | - | - | 407.8 | 0.97 | 321.27 | - | 0.79 |
| DC-Rm-66S-C8 | $f_{c,mortar}$~4 MPa | | C8 | 66x66S | 702.1 | 3.38 | 627.5 | 1.7 | 0.89 |
| DC-C-U | C: Cobblestone, | 400 | - | - | 118.2 | 0.41 | 78.11 | - | 0.66 |
| DC-C-66S-C6 | $f_{c,mortar}$~3 MPa | | C6 | 66x66S | 375.5 | 1.38 | 310.81 | 3.2 | 0.83 |
| DC-C2-U | C2: Cobblestone, | 400 | - | - | 48.20 | 0.46 | 28.82 | - | 0.60 |
| DC-C2-66S-C5 | $f_{c,mortar}$~1.5 MPa | | C5 | 66x66S | 218.3 | 0.52 | 175.64 | 4.5 | 0.80 |

* Mean values of compressive strength, tensile strength and Young modulus, evaluated experimentally, are: 6.3 MPa / 1.1 MPa / 14.4 GPa for C8, 6.7 MPa / 0.8 MPa / 14.4 GPa for C6, 10.1 MPa / 1.4 MPa / 14.4 GPa for C10, 25.7 MPa / 3.0 MPa / 23.4 GPa for C25 and 4.2 MPa / 0.6 MPa / 14.5 GPa for C5.

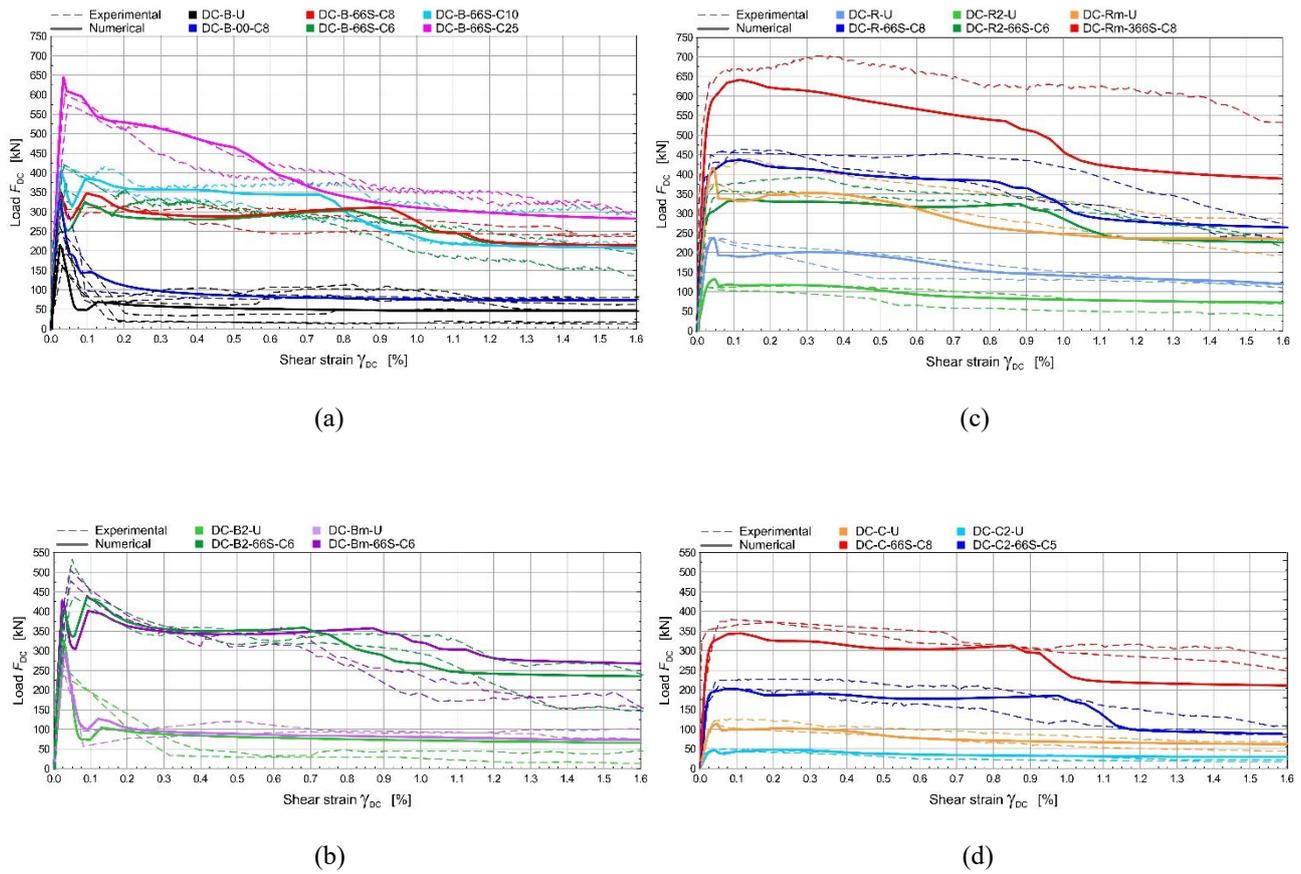

**Fig. 3** Diagonal-compression tests: load-shear strain curves for (a-b) solid brick samples, (c) rubble stones and (d) cobblestone samples. The experimental curves are drawn in dashed lines, the numerical ones with solid lines.

### 3.2. In-plane bending

In-plane, three-point bending tests allow to investigate on the influence of CRM on the in-plane bending mechanism of masonry [13]. The available results concerned 250 mm thick solid brick masonry wallets, with dimensions 780x380 mm$^2$, arranged horizontally and subjected to vertical loading-unloading cycles (span 680 mm) under a constant horizontal pre-compression - Fig. 4a. The tests investigated on the influence of different pre-compression levels and GFRP grid spacing and orientation and also unreinforced samples were tested for comparison. The masonry was type "B" and the mortar of the coating type "C6" (the material characteristics were those already indicated in Table 2). The main features and results are resumed in Table 3 (the samples id. starts with the indication of the type of test, IB, and masonry, B, followed by the GFRP mesh orientation and type and the pre-compression stress level). The backbone curves are also reported in Fig. 4b in terms of applied load, $F_{IB}$, at the varying of the net deflection $\delta_{IB}$ at the intrados.

In the two unreinforced samples (axial stress level 0.15 MPa and 0.30 MPa) a single vertical crack formed abruptly in the mortar joint at the mid-span (Fig. 4c); then, the rocking kinematic mechanism activated, due to the presence of the axial load. The load-deflection curves had an approximately elastic-plastic trend; the higher was the pre-compression, the greater was the resistance $F_{IBu}$. In the CRM strengthened samples, once the vertical crack of the mortar coating occurred at the mid-span (for a cracking load $F_{IBcr}$), the global stiffness decreased but the load continued to increase due to presence of the GFRP mesh, which opposed to the activation of rocking. Other cracks progressively formed in the vicinity, following a vertical trend at the intrados, which tended to incline in the upper portion; a greater crack diffusion was generally encountered for higher reinforcement ratio (Fig. 4d-e). Once the GFRP horizontal yarns collapsed in tension, starting from the lower ones, an abrupt load drop-down occurred. By analyzing the results, it emerged that the resistance increased with the axial stress level (under constant reinforcement ratio) and with the reinforcement ratio (under constant axial stress level). The CRM effectiveness, in terms of resistance, decreased with the axial stress level.

**Table 3**. Experimental in-plane bending tests ($F_{IBcr}$ and $F_{IBu}$ are the first cracking and ultimate load, with respective deflection values, $\delta_{IBcr}$ and $\delta_{IBu}$).

| ID | GFRP mesh Orient.* | Type | Pre-comp. [MPa] | $F_{IBcr}$ [kN] | $\delta_{IBcr}$ [mm] | $F_{IBu}$ [kN] | $\delta_{IBu}$ [mm] | $F_{IBmax(R)}/F_{IBmax(U)}$ |
|---|---|---|---|---|---|---|---|---|
| IB-B-P66S-00 | P | 66x66S | 0.00 | 23.4 | 0.03 | 49.0 | 4.0 | - |
| IB-B-U-15 | - | - | 0.15 | 14.9 | 0.04 | 17.9 | - | - |
| IB-B-P66S-15 | P | 66x66S | | 27.0 | 0.05 | 67.5 | 4.8 | 3.8 |
| IB-B-T99S-15 | T | 99x99S | | 29.4 | 0.06 | 57.6 | 2.9 | 3.2 |
| IB-B-T33S-15 | T | 33x33S | | 33.8 | 0.07 | 99.6 | 4.5 | 5.6 |
| IB-B-U-30 | - | - | 0.30 | 22.2 | 0.05 | 34.2 | - | - |
| IB-B-P66S-30 | P | 66x66S | | 43.2 | 0.08 | 85.1 | 4.9 | 2.5 |
| IB-B-T33S-30 | T | 33x33S | | 44.4 | 0.08 | 108.6 | 5.1 | 3.2 |

* P, T = Parallel/Twisted fiber yarns oriented in the horizontal direction

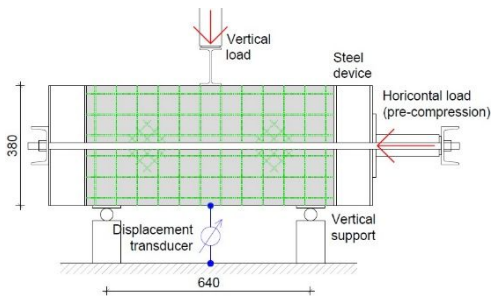
(a)

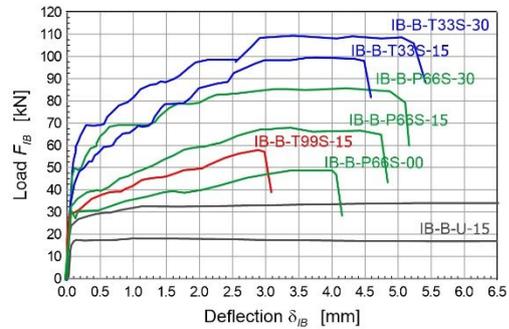
(b)

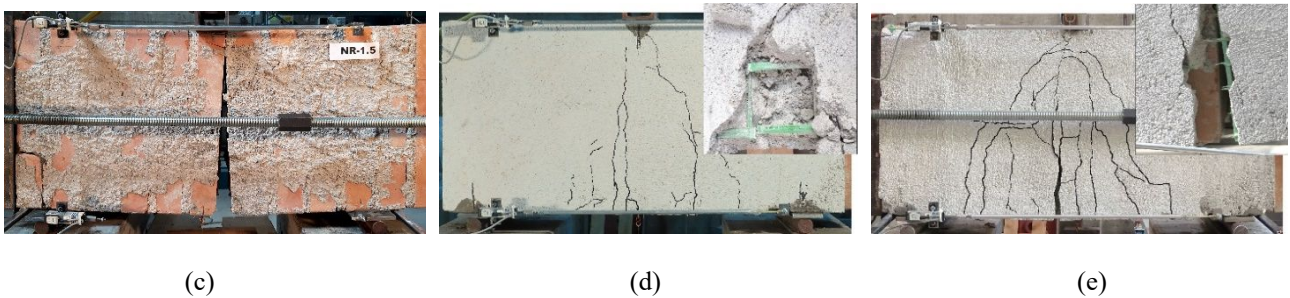

|  (c)  |  (d)  |  (e)  |

**Fig. 4** Experimental in-plane bending tests: (a) test setup, (b) load-deflection backbone curves and crack pattern of the unstrengthened sample (c) IB-B-U-15 and of strengthened samples (d) IB-B-P66S-15 and (f) IB-B-T33S-15, with detail of the lower horizontal yarns at collapse.

### 3.3. Out-of-plane bending

The effectiveness of CRM against the out-of-plane bending mechanism was assessed by means of four-point bending tests [14] performed on masonry samples 3000 mm tall and 1000 mm wide, lying on a vertical linear support, at the mid-thickness, and horizontal linear constraints at the top and at the bottom of the front side (Fig. 5a). Two equal, horizontal forces were applied at the thirds of the height, at the back side. The tests investigated on different masonry types and thickness (with reference to the materials characteristics reported in Table 2, solid brick masonry type "B", 250 mm thick, rubble stone masonry "R", 400 mm thick, and cobblestone masonry "C", 400 mm thick). The mortar of the coating was almost equivalent to type "C8". Unreinforced masonry samples were also tested, as reference. The main test features and results are summarized in Table 4 (the samples id. starts with the indication of the type of test, OB, and masonry, followed by the GFRP mesh orientation and type and by the mortar type). The capacity curves are reported in Fig. 5b in terms of applied load, $F_{OB}$, at the varying of the net deflection $\delta_{OB}$.

Unreinforced samples reached the peak load just before the opening of a horizontal crack in the mid-third of the height (Fig. 5c), at the front side; it followed a sudden resistance reduction and a second phase governed by rocking, due to the self-weight. Thus, the residual resistance resulted higher in thicker samples. Differently, in CRM strengthened samples, once the first cracking appeared horizontally at the front side, the load had a partial reduction but then was gradually recovered and then even increased, as the GFRP vertical yarns crossing the crack progressively carried the tensile stresses. Other parallel cracks progressively formed, almost in the mid-third of the height (Fig. 5d), till the attainment of collapse, due to the yarns breakage. It is generally observed that the sample resistance increased with the sample thickness, while ultimate deflection decreased. Unexpectedly, the cobblestones strengthened sample provided greater resistance than the rubble stone one, despite the same thickness and reinforcement ratio. This can reasonably be attributed to some additional, uncontrolled frictional phenomena that originated in the contact areas with the steel apparatus (e.g. not perfectly clean steel), which also limited the deflection capacity.

**Table 4**. Experimental out-of-plane bending tests ($F_{OBcr}$ and $F_{OBu}$ are the first cracking and ultimate load, with respective deflection values, $\delta_{OBcr}$ and $\delta_{OBu}$).

| ID | GFRP mesh Orient.* | Type | $F_{OBcr}$ [kN] | $\delta_{OBu}$ [mm] | $F_{OBu}$ [kN] | $\delta_{OBu}$ [mm] | $F_{OBu(R)}$ / $F_{OBu(U)}$ |
|---|---|---|---|---|---|---|---|
| OB-B-U | - | - | 9.6 | 0.38 | 3.4** | - | - |
| OB-B-T66S-C8 | T | 66x66S | 36.0 | 32.6 | 45.5 | 32.6 | 13.4 |
| OB-R-U | - | - | 25.5 | 0.46 | 11.6** | - | - |
| OB-R-T66S-C8 | T | 66x66S | 53.1 | 19.9 | 86.2 | 20.2 | 9.0 |
| OB-C-U | - | - | 15.1 | 0.56 | 8.1** | - | |
| OB-C-T66S-C8 | T | 66x66S | 66.4 | 13.8 | 101.3 | 13.8 | 12.5 |

* T = Twisted fiber yarns oriented in the vertical direction.

** Residual load related to the rocking mechanism.

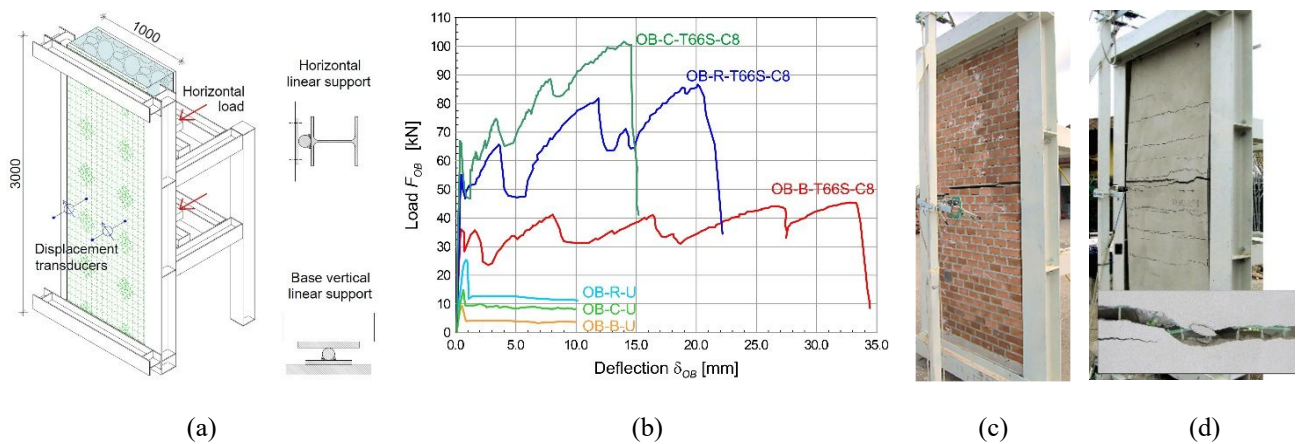

(a)      (b)      (c)      (d)

**Fig. 5** Experimental out-of-plane bending tests: (a) test setup, (b) load-deflection curves and crack pattern of unstrengthened sample (c) OB-B-U and of strengthened samples (d) OB-B-T66S, with detail of the vertical GFRP yarns at collapse in correspondance of the main crack.

## 4. Characteristics of the detailed level model

The numerical simulations on CRM strengthened masonry samples were implemented by using the free, open source, finite element code OOFEM [15]. Nonlinear-static analyses at displacement control were performed considering the nonlinearities of materials and interfaces (Newton-Rapshon solver, with relative displacement and force convergence norms set to 10$^{-3}$). In the prospect of an open science approach, the finite element code adopted and the input files of the models herein described are available for free consultation and use [16]-[17].

The main features of the detailed level numerical model are schematized in Fig. 6: 16.5x16.5x15 mm$^3$ 8-nodes brick elements (*LSpace*) were considered for both the mortar coating and the masonry, while 16.5 mm truss elements (*Truss3d*) for the GFRP yarns. An elastic-brittle behavior in tension was assumed for the yarns (OOFEM material model *Idm1*), with parameters already set in [7] on the basis of experimental direct tensile tests on single yarns (Table 5). Differently, a Concrete-Damage Plasticity model [18], accounting for both cracking and crushing (OOFEM material model *Con2dpm*), was set for the mortar (Table 6), in accordance to the results experimental characterization tests performed on mortar cylinders (compression tests and indirect tensile tests). Line-to-line interface elements (*IntElLine1*), oriented in the yarn direction, tied together the truss elements, representing the yarns, and the edge of the solid elements representing the mortar; point-to-point interfaces (*IntElPoint*) governed the interaction between orthogonal yarns. Face-to-face plane interfaces (*IntElSurfTr*) were introduced to connect the mortar coating to the masonry substrate. All the interfaces (Table 7) were governed by a tangential bond-slip law (*bondceb*), whose parameters were previously set on the basis of experimental tests, as detailed specified in [7]. In particular, pull-out tests of GFRP yarns from the mortar coating, for the line-to-line interfaces, shear tests on mesh nodes, for point-to-point interfaces, and shear-bond tests on CRM coupons applied on masonry wallets, for the face-to-face interfaces.

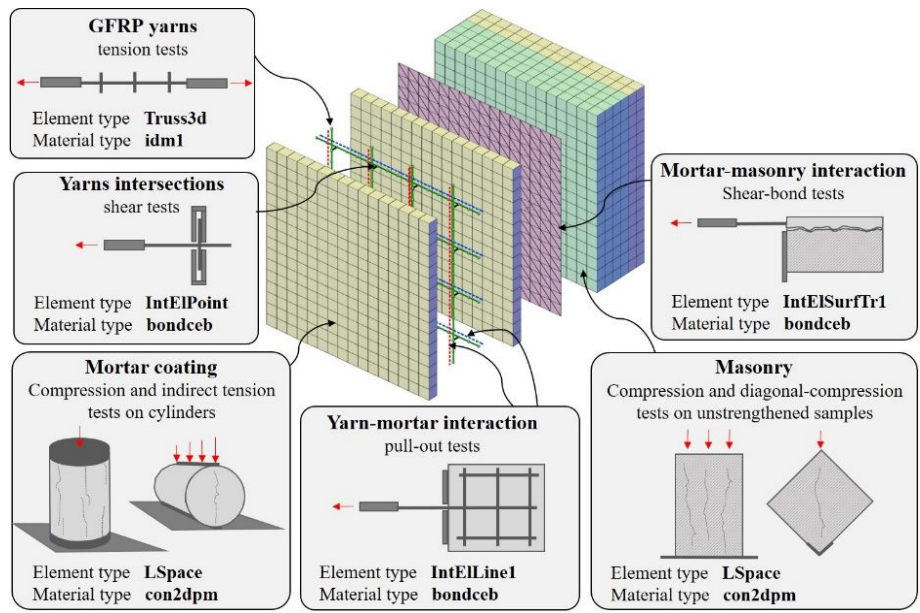

**Fig. 6** Main features of the OOFEM detailed level modelling, with reference to the experimental tests that served for the parameters calibration.

Further details concerning the mentioned elements and materials models are available in the "Manuals" section of the OOFEM webpage [6]. Note that the yarn-mortar interaction properties, previously calibrated by considering experimental pull-out tests with mortar type C8, were herein assumed also for the other mortar types, in lack of related experimental results. Consistently, some authors investigating on TRM characterization evidenced that the yarns bond behaviour mostly depends on their roughness and coating treatment, rather than the mortar matrix characteristics [19]-[21]. Referring to the mortar-masonry interaction, the properties previously calibrated for solid brick masonry were assumed also for stone masonry, lacking of dedicated experimental outcomes; however, this is, reasonably, on the safety side, considering the greater roughness of the stone masonry surfaces.

**Table 5**. Main parameters adopted for the GFRP yarns.

| GFRP yarns (type S) | | | Twisted | Parallel |
|---|---|---|---|---|
| Oofem element type | | | Truss3d | Truss3d |
| Oofem material type | | | Idm1 | Idm1 |
| Cross section | $A$ | [mm$^2$] | 3.8 | 3.8 |
| Young modulus | $E$ | [GPa] | 62.9 | 69.5 |
| Poisson modulus | $n$ | - | 0.3 | 0.3 |
| Softening law | - | | Linear | Linear |
| Peak strain | $\varepsilon_0$ | [%] | 1.88 | 2.12 |
| Strain at 0 stress | $\varepsilon_f$ | [%] | 1.90 | 2.20 |

**Table 6**. Main parameters adopted for the mortar coating.

| Mortar coating | | | C8 | C6 | C10 | C25 | C5 |
|---|---|---|---|---|---|---|---|
| Oofem elem. type | | | LSpace | LSpace | LSpace | LSpace | LSpace |
| Oofem mat. type | | | con2dpm | con2dpm | con2dpm | con2dpm | con2dpm |
| Young modulus | $E$ | [GPa] | 14.4 | 14.4 | 14.4 | 23.4 | 14.5 |
| Poisson modulus | $n$ | [-] | 0.25 | 0.25 | 0.25 | 0.25 | 0.25 |
| Self-weight | $\gamma$ | [kN/m$^3$] | 20.0 | 20.0 | 20.0 | 20.0 | 20.0 |
| Comp. strength | $f_c$ | [MPa] | **6.29** | **6.71** | **10.14** | **25.74** | **4.23** |
| Tens. strength | $f_t$ | [MPa] | **1.10** | **0.8** | **1.36** | **2.97** | **0.55** |
| Dilation | $\psi$ | [°] | 40 | 40 | 40 | 40 | 40 |
| Softening law | - | [-] | Linear | Linear | Linear | Linear | Linear |
| Hardening param. | $b_h$ | [-] | 0.002 | 0.002 | 0.002 | 0.002 | 0.004 |
| | $h_p$ | [-] | 0.0 | 0.0 | 0.0 | 0.0 | 0.0 |
| Softening param | $w_f/h$ | [-] | 0.011 | 0.013 | 0.009 | 0.003 | 0.012 |
| | $a_{soft}$ | [-] | 4.0 | 4.0 | 10.0 | 15.0 | 4.0 |

**Table 7**. Main parameters adopted for the interfaces.

|  |  | Yarn-mortar interaction | | Yarns intersections | | Masonry-mortar interaction |
|---|---|---|---|---|---|---|
|  |  | Twisted | Parallel | Twisted | Parallel |  |
| Oofem element type | | intElLine1 | intElLine1 | intElPoint | intElPoint | IntElSurfTr1 |
| Oofem material type | | bondceb | bondceb | bondceb | bondceb | bondceb |
| Equiv. perimeter t | [mm] | 9.57 | 18.00 | - | - | |
| Normal stiffness $k_n$ | [N/mm$^2$] | $10^3$ | $10^3$ | 0 | 0 | $10^4$ |
| Tang. stiffness $k_t$ | [N/mm$^2$] | $10^3$ | $10^3$ | $10^4$ N/mm | $10^4$ N/mm | 500 |
| Exponent $\alpha$ | [-] | 0.4 | 0.4 | 0.4 | 0.4 | 0.55 |
| Peak bond $\tau_{max}$ | [MPa] | 3.30 | 2.00 | 458.0 | 550.0 | 1.18 |
| Residual bond $\tau_f$ | [MPa] | 2.20 | 0.25 | 0.0 | 0.0 | 0.02 |
| Slip parameters $s_1$ | [mm] | 0.1 | 0.1 | 0.5 | 0.5 | 0.04 |
| $s_2$ | [mm] | 0.1 | 0.1 | 1.5 | 0.5 | 0.05 |
| $s_3$ | [mm] | 1.2 | 1.0 | 10.0 | 0.6 | 0.1 |

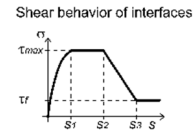

Shear behavior of interfaces

The masonry was considered as an equivalent homogeneous material, according to a macro-modelling approach [22]; its behaviour was simplistically assumed isotropic and governed by the Concrete-Damage Plasticity model [18] (OOFEM material model *Con2dpm*). Six masonry types were analysed, both regular (B, B2, R) and irregular (R2, C, C2); the parameters, summarized in Table 8, were calibrated on the basis of the results obtained by testing unstrengthened masonry elements in diagonal compression [10]-[12] and in compression [14].

In particular, the masonry tensile strength $f_t$ was calculated from the results of diagonal compression tests (see subsection 3.1), by applying the formulation $f_t = \alpha\, F_{DCmax} / bt$, being $F_{DCmax}$ the peak diagonal load (mean value obtained from the experimental tests), and $bt$ the sample cross section. Accordingly to Brignola et al. [23], $\alpha$ was ranged from 0.5 (for regular masonry, characterized by a brittle behaviour) to 0.35 (for irregular masonry, where a relevant redistribution of the stresses due to the interlocking effect occurs). The softening parameters in tension were set so to fit the mean experimental results of the experimental diagonal compression tests performed on unstrengthened samples. For solid brick masonry, the $w_f/h$ parameter, governing the behavior in tension, was set quite low (Table 8), so to correctly catch the abrupt drop down of resistance (brittle failure). Differently, stone samples necessitated a softer tensile degradation (higher values of the $w_f/h$ parameter), accounting for a more pronounced interlocking effect among blocks. It is worth to note that, since obtained from diagonal compression tests, the derived parameters are referred to a loading direction angle of 45° with respect to the bed joints; but, due to the hypothesis of homogeneous isotropic material, they are assumed constant in all directions. Although this simplified hypothesis would not be exactly adequate to accurately simulate the behavior of unreinforced masonry, it is considered an acceptable approximation for the purposes of this study, focused on reinforced masonry, whose behavior is mainly governed by the CRM system.

The masonry Young modulus and compression parameters were set so to fit the mean experimental results of experimental compression tests on plain masonry wallets available in the literature. But, for regular masonry, the compressive strength is strongly influenced by the loading direction, as proved by Page [24]; therefore, since a homogeneous isotropic behaviour was assumed in the numerical models, averaged values were considered (consistently with Page's results). For the purpose of this study, focusing on strengthened masonry, these simplified assumptions were acceptable.

**Table 8.** Main parameters adopted for the masonry.

| MASONRY | | Solid brick | | Rubble stone | | Cobblestone | |
|---|---|---|---|---|---|---|---|
| | | **B** | **B2** | **R** | **R2** | **C** | **C2** |
| Oofem elem. type | | LSpace | LSpace | LSpace | LSpace | LSpace | LSpace |
| Oofem mat. type | | con2dpm | con2dpm | con2dpm | con2dpm | con2dpm | con2dpm |
| Young modulus $E$ | [GPa] | 4.27 | 5.30 | 2.43 | 1.42 | 1.26 | 0.60 |
| Poisson modulus $N$ | [-] | 0.45 | 0.45 | 0.45 | 0.45 | 0.45 | 0.45 |
| Self-weight $\gamma$ | [kN/m³] | 18.0 | 18.0 | 21.0 | 21.0 | 19.0 | 19.0 |
| Comp. strength $f_c$ | [MPa] | **5.12** | **7.00** | **2.13** | **1.20** | **1.04** | **0.45** |
| Tens.strength $f_t$ | [MPa] | **0.320** | **0.440** | **0.208** | **0.104** | **0.089** | **0.036** |
| Dilation $\Psi$ | [°] | 30 | 30 | 35 | 40 | 40 | 40 |
| Softening law - | [-] | Linear | Linear | Linear | Linear | Linear | Linear |
| Hardening param. $b_h$ | [-] | 0.003 | 0.003 | 0.006 | 0.006 | 0.006 | 0.006 |
| $h_p$ | [-] | 0.0 | 0.0 | 0.0 | 0.0 | 0.0 | 0.0 |
| Softening param $w_f/h$ | [-] | 0.0001 | 0.0002 | 0.004 | 0.003 | 0.004 | 0.003 |
| $a_{soft}$ | [-] | 5.0 | 8.0 | 15.0 | 15.0 | 10.0 | 10.0 |

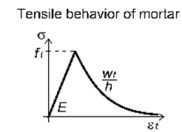

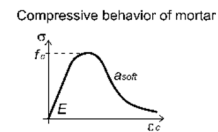

## 5. Numerical results

The detailed level numerical model was applied to the simulations of CRM strengthened masonry elements subjected to diagonal compression (DC), in-plane bending (IB) and out-of-plane bending (OB) - Fig. 7, according to the characteristics and setup of experimental tests previously described in section 0. To reduce the computational effort, symmetry planes were considered, so to allow the modelling only of some portions of the specimens. Translational constrains were applied to the nodes owning to the symmetry planes, avoiding displacements orthogonal to such planes. Moreover, the walls swelling was avoided, since no evident through-the-thickness deformations or damage emerged in all the considered experimental tests.

The total number of elements/interfaces in the DC and OB strengthened models is about 18-19 thousands and 7.5 thousands in the IB models. To run an analysis, by using 10 cores parallel processors Intel® Xeon® CPU E5-2630 v3 @ 2.40 GHz, takes about 4 hours for DC, 10 hours for IB and 24 hours for OB.

The equivalent thickness set for the CRM layer, to be intended as the actual average thickness, was 30 mm for solid brick masonry (B and B2), 40 mm for regular rubble stone masonry (R) and 50 mm for irregular rubble stone and cobblestone masonry (R2, C and C2), consistently with [11] and [14]. These differences, actually found in the samples, are related to the different irregularity of the masonry surface: a higher roughness implies also a greater thickness variability in respect to the minimum nominal value, evaluated from the most protruding stones.

Based on the experimental outcomes available for comparison, the simulations of diagonal compression tests focused on the evaluation of the reliability of the numerical model at the varying of the masonry type and thickness and of the mortar type, while those concerning in-plane bending tests, on the influence of the reinforcement ratio and of the pre-compression level. Those related to out-of-plane bending investigated the influence of the masonry type and thickness.

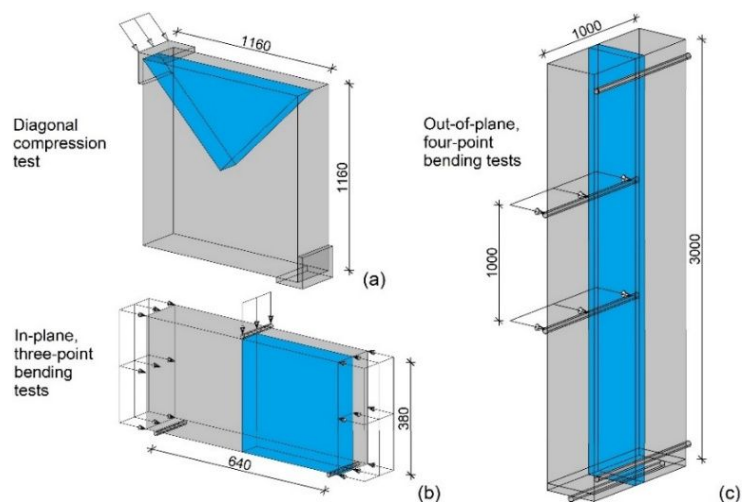

**Fig. 7** Schematization of (a) diagonal compression, (b) in-plane bending and (c) out-of-plane bending tests, highlighting the modelled portions of the samples.

5.1. Diagonal compression tests

In the simulations of the diagonal compression tests, one-eighth of the sample was modelled (Fig. 7a) and the diagonal displacement was applied in correspondance of the corner bracket. The load - shear strain $F_{DC}$-$\gamma_{DC}$ numerical curves are plotted with solid lines in Fig. 3; the main results are summarized in Table 9.

The initial linear behavior of the unreinforced masonry samples was interrupted by the tensile failure at the center of the panel, originating a diagonal crack and, then, the formation of a pushing masonry wedge around the corner device (Fig. 8a). For solid brick masonry, an abrupt drop down of resistance occurred. Differently, in stone samples, the resistance degradation resulted smoother, as a consequence of the adopted softening law in tension, combined with the compressive resistance of the diagonal masonry struts.

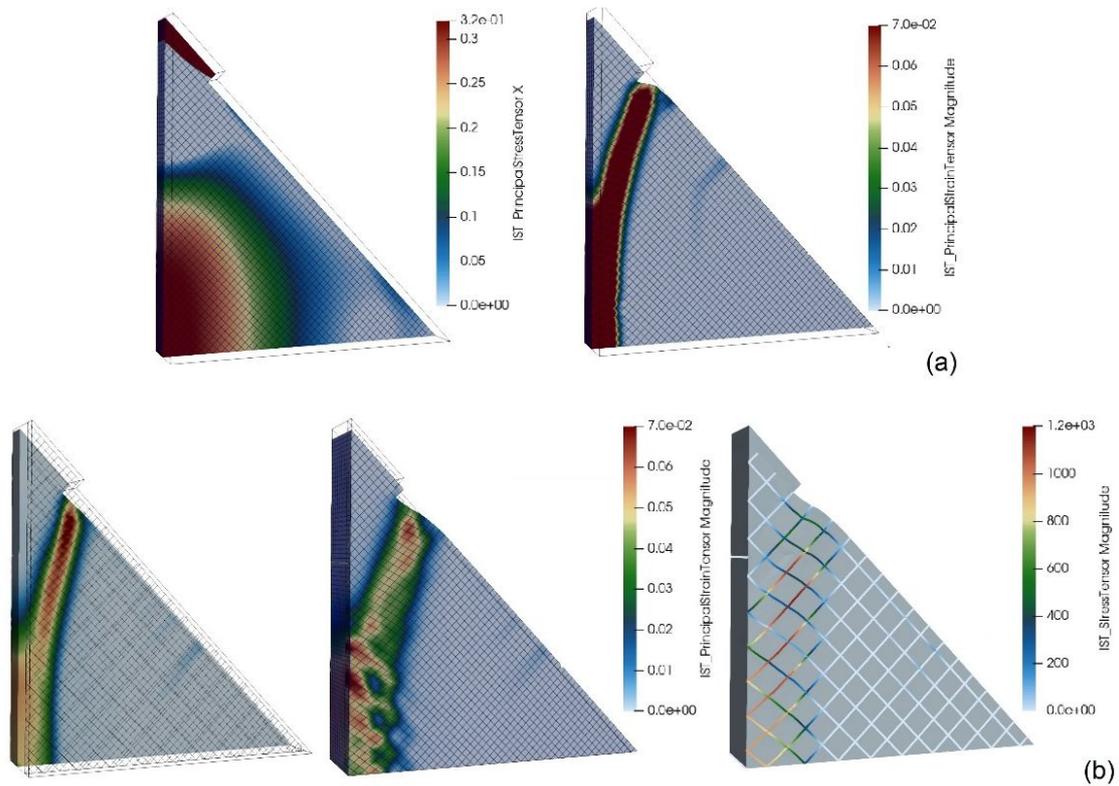

**Fig. 8** Diagonal compression tests: (a) principal tensile stresses at peak load and crack pattern at $\gamma_{DC}$ = 0.8% in unstrengthened sample DC-B-U and (b) crack pattern in masonry and in mortar coating and tensile stresses in the GFRP wires at $\gamma_{DC}$ = 0.8% in strengthened sample DC-B66S-C8.

In the CRM strengthened sample, significant performance improvements emerged in respect to plain masonry. The peak load was detected as the tensile resistance in the mortar coating and in the masonry were reached, almost simultaneously, in the center of the panel. At the increasing of the applied diagonal displacement, the damage gradually spread; actually, in respect to the masonry, the damaged area in the mortar coating appeared wider, but with lower tensile strains, due to the diffusive action of the embedded GFRP mesh (Fig. 8b). The failure of the GFRP yarns started at about $\gamma_{DC}$ = 0.8-0.9% and determined an evident drop of resistance. Focusing on DC-B-66S tests (Fig. 3a), it is observed that the peak load increased with the tensile strength of the mortar coating but the load decrease was sharper. In fact, the first part of the tests, before the peak, was mostly influenced by the matrix characteristics, while the GFRP mesh plays a negligible role; but then, in the post-peak phase, the global behavior was governed by the GFRP mesh and its capability in opposing to the cracks opening. This was clearly noticed by comparing DC-B-66S-C8 with DC-B-00-C8 (sample without GFRP mesh), which attained almost the same peak load but incurred to pseudo-ductile and brittle collapse mode, respectively. By comparing the results obtained from the different solid brick masonry types with same CRM characteristics (mesh 66S and mortar C6) - Fig. 3a-b , it is observed that the effectiveness of the CRM in terms of resistance tended to reduce at the increasing of both the masonry properties (B2 vs. B) and thickness (Bm vs. B). The shear strain associated to the yarns failure was not particularly affected by the masonry thickness (Bm vs. B), while seems to emerge that a stronger masonry (B2 vs. B) or mortar (C25 vs. C10 or C8) anticipates this occurrence. These outcomes were confirmed also from the comparison of results on the stone masonry samples DC-R-66S-C8, DC-Rm-66S-C8, DC-R2-66S-C6 (Fig. 3c) and the cobblestones sample DC-C-66S-C8 (Fig. 3d).

**Table 9**. Main numerical results of diagonal-compression tests ($F_{DCmax}$ is the peak load and $\gamma(F_{DCmax})$ the shear strain at peak load, $F_{DC8}$ the load in correspondance of a shear strain $\gamma_{DC} = 0.8\%$).

| ID | $F_{DCmax}$ [kN] | $\gamma(F_{DCmax})$ [⁰/₀₀] | $F_{DC8}$ [kN] | $F_{DCmax(R)}/F_{DCmax(U)}$ | $F_{DC8}/F_{DCmax}$ |
|---|---|---|---|---|---|
| DC-B-U | 214.96 | 0.25 | 48.3 | - | 0.22 |
| DC-B-00-C8 | 330.40 | 0.26 | 78.3 | 1.5 | 0.24 |
| DC-B-66S-C8 | 350.48 | 0.25 | 308.3 | 1.6 | 0.88 |
| DC-B-66S-C6 | 323.01 | 0.23 | 304.6 | 1.5 | 0.94 |
| DC-B-66S-C10 | 404.68 | 0.28 | 292.6 | 1.9 | 0.72 |
| DC-B-66S-C25 | 645.12 | 0.34 | 334.0 | 3.0 | 0.52 |
| DC-B2-U | 293.51 | 0.27 | 83.0 | - | 0.28 |
| DC-B2-66S-C6 | 436.92 | 0.25 | 305.3 | 1.5 | 0.70 |
| DC-Bm-U | 326.73 | 0.26 | 74.56 | - | 0.23 |
| DC-Bm-66S-C6 | 428.12 | 0.23 | 354.15 | 1.3 | 0.83 |
| DC-R-U | 237.80 | 0.46 | 151.74 | - | 0.64 |
| DC-R-66S-C8 | 437.16 | 1.15 | 384.74 | 1.8 | 0.88 |
| DC-R2-U | 132.55 | 0.49 | 83.8 | - | 0.63 |
| DC-R2-66S-C6 | 335.61 | 1.18 | 321.0 | 2.5 | 0.96 |
| DC-Rm-U | 416.16 | 0.46 | 265.8 | - | 0.64 |
| DC-Rm-66S-C8 | 641.56 | 1.18 | 539.1 | 1.5 | 0.84 |
| DC-C-U | 113.76 | 0.47 | 71.0 | - | 0.62 |
| DC-C-66S-C6 | 344.31 | 1.11 | 311.1 | 3.0 | 0.90 |
| DC-C2-U | 47.71 | 0.39 | 33.2 | - | 0.70 |
| DC-C2-66S-C5 | 272.75 | 0.70 | 250.6 | 5.7 | 0.92 |

### 5.2. In-plane bending tests

In the numerical model developed for the simulation of in-plane, three-point bending tests, just one-fourth of the sample was modelled, accounting for the sample symmetry (Fig. 7b). The vertical displacements were constrained at the support nodes; horizontal forces were at first applied at the lateral free end and maintained constant, to represent the pre-compression load; then, the vertical displacement was applied at the top, at the mid span. The load deflection $F_{IB}$-$\delta_{IB}$ numerical curves are plotted in Fig. 9 with solid lines; the main results are also reported in Table 10.

The simulations of plain masonry elements exhibited a single crack in the mid span (at the intrados) and were then governed by the rocking kinematic mechanism (which activated only in presence of a pre-compressive load), resulting in a plastic branch of the $F_{IB}$-$\delta_{IB}$ curves. In the simulations of CRM strengthened samples, a vertical crack firstly appeared at the mid span, at the intrados, due to bending, as the masonry and the mortar of the coating attained almost simultaneously to their respective tensile strengths. Then, as the deflection increased, other cracks originated from the intrados, in 33S and 66S configurations; their slope tended to vary in the upper portion of the sample, due to shear stresses (Fig. 10a-b); differently, in 99S configurations, just one small additional crack formed after the first cracking (Fig. 10c). The maximum load was attained when the ultimate deformation was reached in the horizontal yarn at the intrados, at the mid span. Strengthened samples performed almost elastic-plastic curves with hardening and attained to higher resistance values in respect to unreinforced ones, due to the tensile contribution of the GFRP yarns, which opposed to the free rocking until breakage. The ultimate deflection was not significantly influenced by the stress level and was also comparable in 66S and 33S samples; differently, 99S samples evidenced significantly lower displacement capacities. In fact, looking at the relative yarn-mortar displacements in correspondance of the lower horizontal yarn, at peak load, the slips resulted very low for 33S (~0.1 mm) and 66S (~0.2 mm) configurations, while reached higher values (~1.6 mm) for 99S. Moreover, consistently with the different crack patterns, as the grid pitch widened, the yarn length affected by high stresses reduced.

Summarizing, the numerical simulation of in-plane bending tests confirmed that: 1) with the same reinforcement ratio, both the first cracking and the maximum load increase as the axial stress increases; 2) with the same axial stress level, the maximum load increases with the reinforcement ratio; 3) the reinforcement effectiveness increases by decreasing the axial stress; 4) a closer grid pitch promotes the cracks diffusion and can prevent excessive slips in yarns, that anticipates failures.

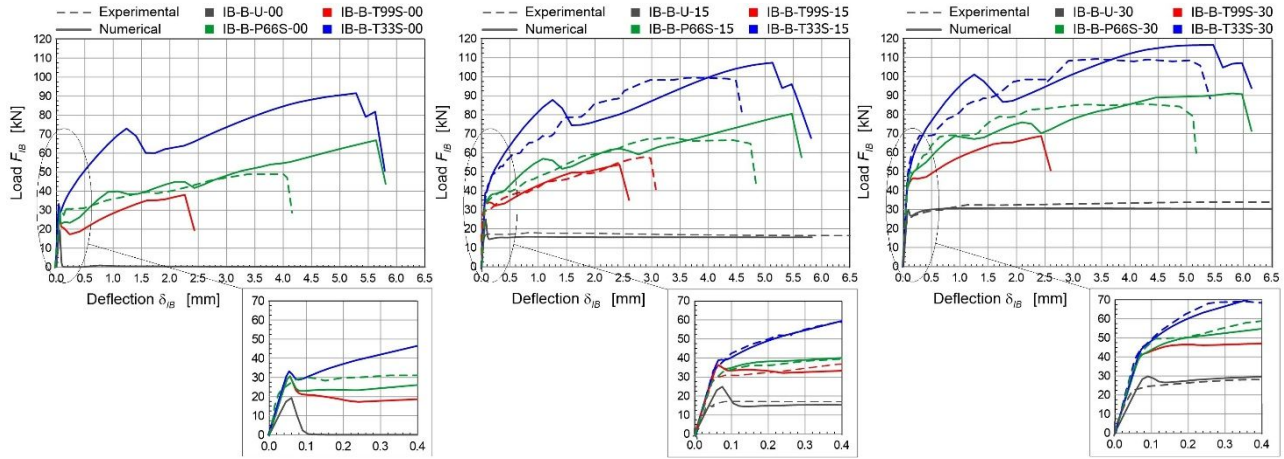

**Fig. 9** In-plane bending capacity curves: comparison between numerical (solid lines) and experimental (dashed lines) results.

**Table 10**. Main numerical results of in-plane bending tests ($F_{IBcr}$ and $F_{IBu}$ are the cracking and ultimate load, with respective deflections, $\delta_{IBcr}$ and $\delta_{IBu}$).

| ID | $F_{IBcr}$ [kN] | $\delta_{IBcr}$ [mm] | $F_{IBu}$ [kN] | $\delta_{IBu}$ [mm] | $F_{IBmax(R)}/F_{IBmax(U)}$ |
|---|---|---|---|---|---|
| IB-B-U-00 | 19.32 | 0.06 | 0.00 | - | - |
| IB-B-P66S-00 | 30.6 | 0.06 | 66.8 | 5.64 | / |
| IB-B-T99S-00 | 30.6 | 0.06 | 38.1 | 2.26 | / |
| IB-B-T33S-00 | 33.2 | 0.06 | 91.6 | 5.30 | / |
| IB-B-U-15 | 24.9 | 0.08 | 15.5 | - | - |
| IB-B-P66S-15 | 36.2 | 0.06 | 80.5 | 5.47 | 5.2 |
| IB-B-T99S-15 | 36.3 | 0.06 | 53.7 | 2.43 | 3.5 |
| IB-B-T33S-15 | 38.8 | 0.06 | 107.4 | 5.47 | 6.9 |
| IB-B-U-30 | 29.8 | 0.09 | 30.5 | - | - |
| IB-B-P66S-30 | 41.8 | 0.08 | 91.1 | 5.81 | 3.0 |
| IB-B-T33S-30 | 44.0 | 0.07 | 116.6 | 5.47 | 3.8 |
| IB-B-T99S-30 | 41.7 | 0.08 | 68.7 | 2.43 | 2.3 |

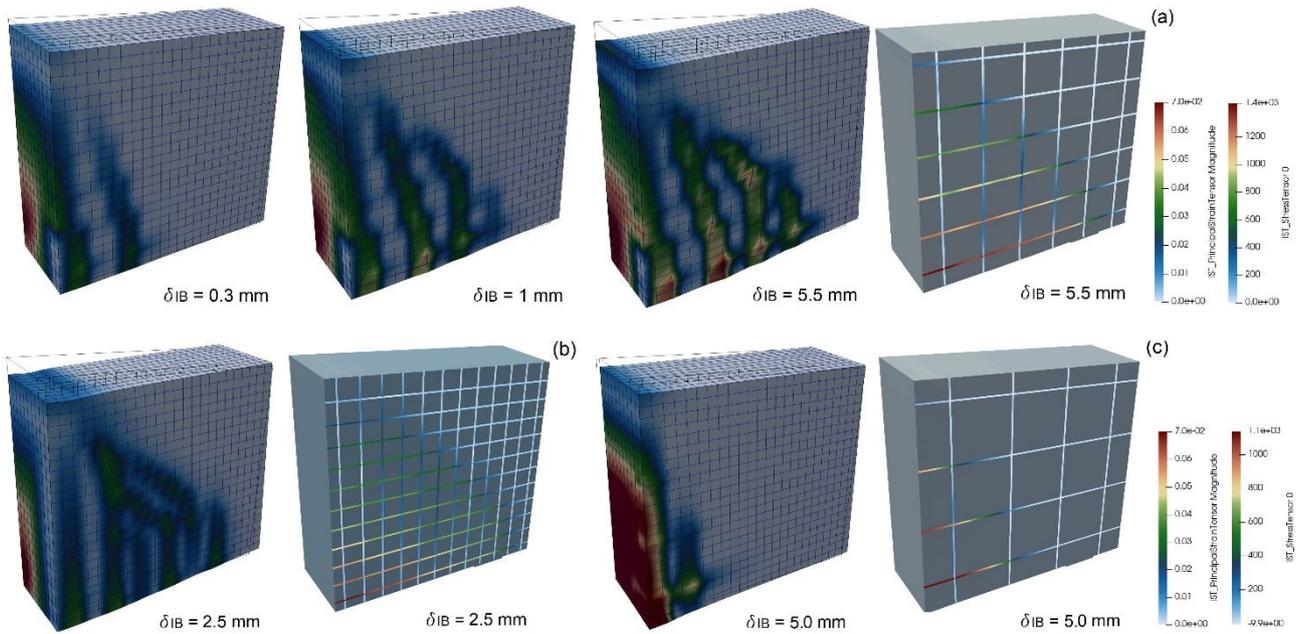

**Fig. 10** In-plane bending numerical simulations: crack pattern and tensile stresses in the GFRP wires at collapse for (a) IB-B-15-P66S, (b) IB-B-15-T33S and (c) IB-B-15-T99S.

5.3. Out-of-plane bending tests

In the simulations of out-of-plane, four-point bending tests, a wall strip of 66 mm width was modelled, since bi-dimensional bending effects in the experimental tests resulted negligible (Fig. 7c). The horizontal displacements were constrained in correspondance of the top and bottom supports, at the front side. Coherently with the experimental setup, a rigid slab was modelled at the bottom and the vertical displacement was avoided at the base, in the mid of the wall thickness. Point-to-point interface elements accounting for the apparatus steel-to-steel friction were introduced at the supports and loading points (coefficient of static friction 0.30, according to [25]). Coherently to the experimental loading procedure, the sample self-weight, in the vertical direction, was at first applied and maintained constant during the simulation; then, horizontal forces were introduced at the thirds of the sample height, at the back side. Direct displacement control was performed by increasing monotonically the wall deflection. The numerical results are reported with solid lines in Fig. 11a in terms of applied load, $F_{OB}$, at the varying of the net deflection at the midspan, $\delta_{OB}$; the main results are summarized in Table 11.

The peak load in the unreinforced samples occurred just before the masonry first cracking, at two-thirds of the wall height ( front side); then, the load rapidly dropped down to a residual load (related to the rocking kinematic mechanism), which resulted higher for the thicker samples. In CRM strengthened samples, the first cracking occurred in the mortar coating at two-thirds of the wall height; then, several horizontal, parallel cracks sequentially originated from the front side (Fig. 11b), mainly in the mid-third of the height. The portions of the vertical yarns crossing the cracks attained to high stresses;

the failure of the yarns in correspondance of the most stressed portion led to the sample collapse; meanwhile, no debonding phenomena emerged at the mortar-masonry interface. Strengthened samples reached significantly higher resistance values in respect to plain masonry, since the GFRP reinforcement opposed to free rocking until breakage. No significant variations occurred in the strengthened samples by varying the masonry type (OB-R-T66S-C8 and OB-C-T66S-C8), indicating that the global behavior is mainly governed by the CRM; a decrease in the masonry thickness (OB-B-T66S-C8) resulted in a lower resistance and a higher deflection. The simulation OB-C-T66S-C8f, for the cobblestones sample, evidenced that a doubled friction coefficient at the apparatus contact points (e.g. due to not perfectly cleaned steel devices) can provide a greater resistance and lower deflection.

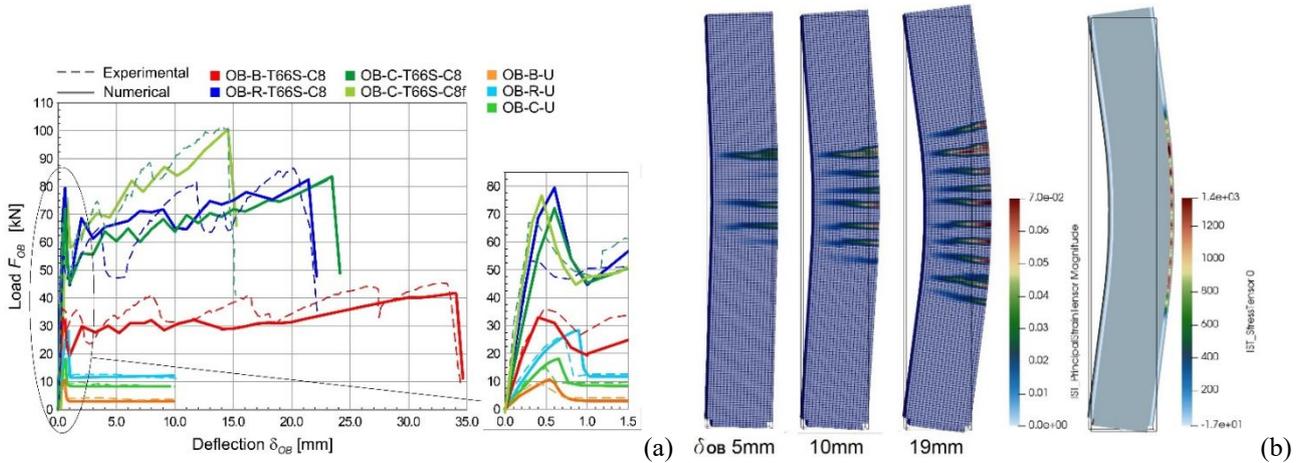

**Fig. 11** Out-of-plane bending tests: (a) capacity curves, with comparison between numerical (solid lines) and experimental results (dashed lines), and (b) example of crack pattern evolution and of tensile stresses in the GFRP yarns just before collapse (sample OB-R-T66-C8).

**Table 11.** Main numerical results of out-of-plane bending tests ($F_{OBcr}$ and $F_{OBu}$ are the first cracking and ultimate load, with respective deflection values, $\delta_{OBcr}$ and $\delta_{OBu}$).

| ID | GFRP mesh Orient.* | GFRP mesh Type | $F_{OBcr}$ [kN] | $\delta_{OBu}$ [mm] | $F_{OBu}$ [kN] | $\delta_{OBu}$ [mm] | $F_{OBu(R)}$ / $F_{OBu(U)}$ |
|---|---|---|---|---|---|---|---|
| OB-B-U | - | - | 10.7 | 0.55 | 3.2 | - | |
| OB-B-T66S-C8 | T | 66x66S | 33.1 | 0.42 | 41.7 | 34.1 | 13.0 |
| OB-R-U | - | - | 28.8 | 0.93 | 12.2 | - | |
| OB-R-T66S-C0 | T | 66x66S | 79.9 | 0.6 | 81.9 | 21.5 | 6.7 |
| OB-C-U | - | - | 18.2 | 0.68 | 8.0 | - | |
| OB-C-T66S-C8 | T | 66x66S | 72.1 | 0.6 | 83.2 | 23.6 | 10.4 |
| OB-C-T66S-C8f* | T | 66x66S | 76.4 | 0.46 | 100.8 | 14.6 | 12.6 |

## 6. Comparison with experimental results

To assess the reliability of the simulations, the numerical curves (drawn in solid lines in Fig. 3, Fig. 9 and Fig. 11a) were compared with the experimental ones (reported with dashed lines). The main representative load and strain (or displacement) values are also compared in the bar charts plotted in Fig. 12, for diagonal compression tests and Fig. 13 for in-plane and out-of-plane bending tests.

In the DC simulations, the error, for both $F_{DCmax}$ and $F_{DC8}$, fell within the range ±15% in about 90% of the CRM strengthened samples. Actually, looking at the capacity curves, the load decrease at the occurrence of the GFRP yarns failure generally resulted less gradual than in the experiments; this was reasonably due to the mean values of the mechanical parameters assumed in the simulations for the yarns; differently, some scatter affects their actual performances. Also the slightly different behavior of the twisted and parallel fibers yarns may have smoothed the load drop in the tests (instead, in the numerical models, the mesh orientation was assumed symmetrical).

For the OB tests, the prediction of the ultimate point fell within the range of a ±10% error, even though the cracking load tended to be overestimated in both the stone samples and the trend of the numerical curves in the post-cracking phase resulted smoother than the experimental ones. The greater thickness variability of the mortar coating, detected in the experimental samples, instead of the constant value assumed in the simulations, is believed to be the main reason.

Acceptable errors resulted also in the IB models, although it was often recognized a general tendency to overestimate the failure point, in terms of both force and deflection. Since the experimental tests were performed after a long period from the sample construction and reinforcement application (almost 10 years), it seems reasonable to hypothesize some degradation in the materials and, in particular in the yarn-mortar bonding, as well as in the yarns resistance. These degradations can justify a reduction of the CRM effectiveness (this aspect will be furtherly investigated through sensitivity analysis in the next section).

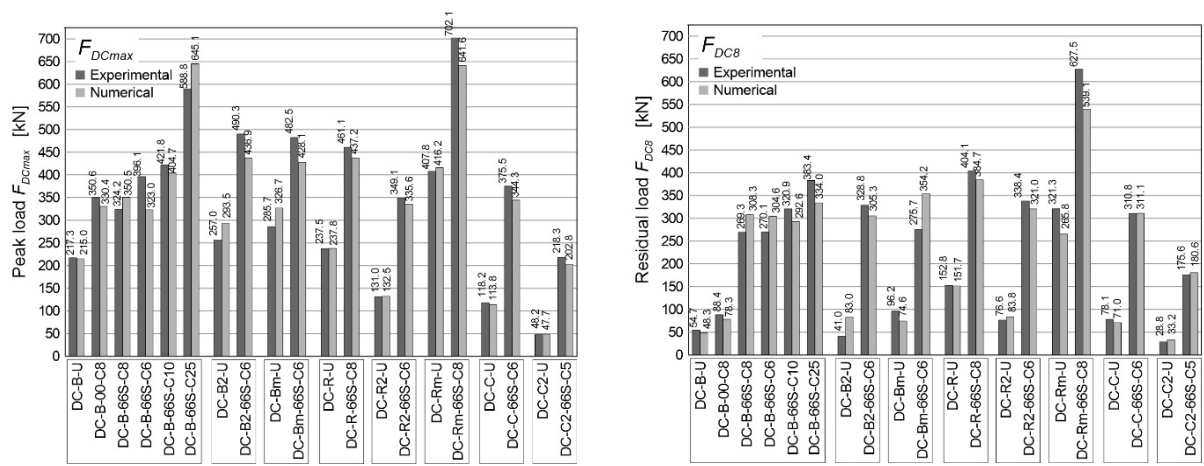

**Fig. 12** Diagonal compression tests: experimental-numerical comparison in terms of peak load and residual load at $\gamma_{DC}$ = 0.8%.

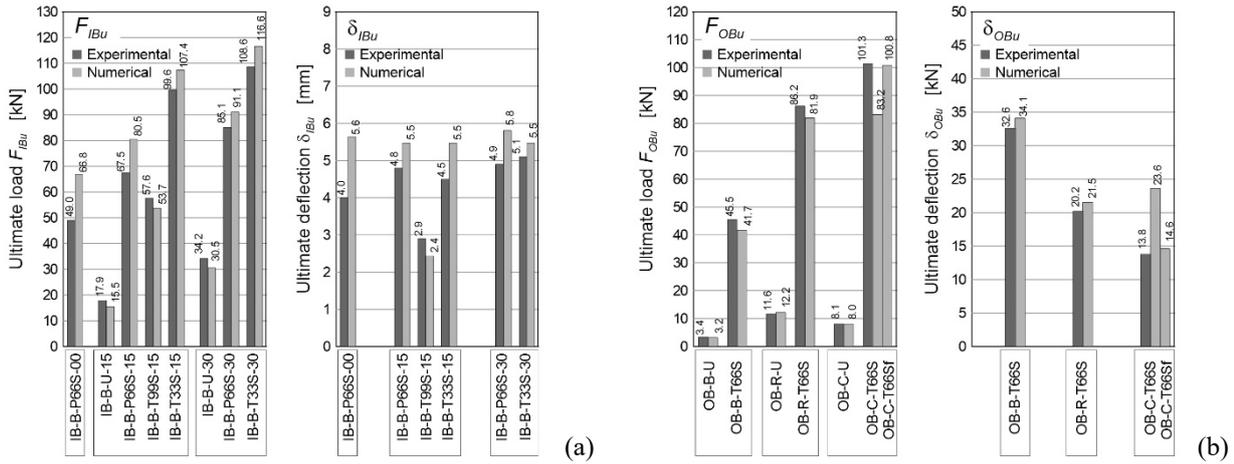

**Fig. 13** Bending tests: experimental-numerical comparison in terms of ultimate load $F_u$ and deflection $\delta_u$, for (a) in-plane loading and (b) out-of-plane loading.

However, globally comparing the experimental (section 3) and numerical outcomes (section 5), it emerged a good accordance in terms of damage evolution and failure modes, proving the capability of the numerical model in detecting the different collapse mechanisms in CRM strengthened masonry, related to both diagonal cracking and bending (in-plane and out-of-plane). Also the capacity curves trends generally look comparable. This is a significant result, since highlights the robustness of the developed numerical model, considering that the material parameters were calibrated on the basis of other experimental characterization tests, rather than suitably set to fit the results obtained from the specific experimental tests on a CRM strengthened masonry element. Moreover, the outputs proved to be reliable at the varying of the characteristics of load pattern, masonry and reinforcement. However, as discussed, some discrepancies emerged. These discrepancies can realistically be attributed to some assumptions adopted in the numerical models. Firstly, it has to be considered that, in the numerical simulation, the nominal, mean values of the materials and interfaces parameters, obtained experimentally from sampling characterization tests, were assumed; but, in the experimental masonry samples, these can actually be affected by significant variations from one sample to another, especially in materials with a brittle behavior. Besides, the masonry intrinsic heterogeneity was neglected in the numerical model (homogeneous material assumption) and the masonry tensile strength was derived from the results of experimental diagonal compression tests on plain masonry and assumed uniform in all directions (isotropic material assumption), even though in the investigated bending configurations (IB, OB) the cracks were parallel to the masonry bed joints. Moreover, the experimental tests concerned loading-unloading cycles, that can imply some cumulative damage with degrading effect on the global performances; differently, the numerical simulations were monotonic and the calibration of the CRM behavior was achieved on the basis of monotonic characterization tests. Actually, the effect of the cyclic loading can lead to a greater deformability (as evidenced e.g. in [26]-[28]). Furthermore, while the numerical simulations were conducted under a precise "displacement control", this was generally not possible in the experimental tests, due to the spring-back effects in the adopted setups, especially when abrupt drops of resistance occurred. Lastly, it should be considered that, for each configuration, one experimental test (IB and OB setups) or at most a couple (DC setup) were performed, due to the onerousness of the full-scale samples; thus, the reference results should be considered as an indicative trend rather than in a purely quantitative terms.

## 7. Sensitivity analysis

Parametric analyses were performed for the different testing setups, to evaluate the sensitivity of the numerical results to variations in the mechanical and geometrical characteristic of the samples. Typically, ±33% variations on the main material parameters were investigated.

For diagonal compression tests, sample DC-B-66S-C8 was selected as reference test; the results are compared in Table 12. Generally, the variations in the GFRP yarns parameters (Fig. 14a) did not affect the initial stiffness and the peak load of the sample, confirming the secondary role of the embedded reinforcement until the cracking of the mortar matrix. When the reduction/increase of the peak and ultimate tensile strains was performed (id. YT), the shear strain corresponding to the beginning of the yarns failure (called "failure point" in the following) reduced/increased without affecting the curve trend up to that point. At the increase/decrease of the yarns Young modulus (YY), the ordinates of the post-peak branch rose/fell and the failure point reached a higher/lower force level and a lower/higher shear stain. The reduction of the dry fibers area in yarns (YA) resulted in a combination of the YT and YY effects: a lowering of the ordinates of the post-peak branch and an anticipation of the failure point, in terms of both force and shear strain; the increase of the same parameter led to higher post-peak ordinates but did not affect significantly the shear strain at failure. By switching the parallel and twisted fiber yarns (YO) the failure point was postponed, since the parallel fiber yarns (having slightly higher stiffness and strength - Table 1) resulted the most stressed in this new configuration.

When varying the mortar coating parameters (Fig. 15a), the increase/decrease of the tensile strength (CT) reflected in higher/lower ordinates, starting from the peak load, under almost the same shear strain levels; small variations were also detected in the failure point. Similar effects were noted for the increase/decrease of the compressive strength (CC). The increase/decrease of the mortar Young modulus (CY) produced a slight increase/decrease the initial stiffness and a lowering/rising of the peak load, but had a limited influence on the failure point.

When referring to the masonry parameters (Fig. 15a), it is observed that the increase of the tensile strength (MT) induced a higher peak load, while did not affect significantly the post-peak behavior. Differently, a decrease of the same parameter resulted in a lower first peak load and in a steeper load decrease, but then the load was recovered and the failure point was rather maintained. The increase/decrease of the compressive strength (MC) did not alter the first peak load but significantly increased/decreased the ordinates of the post peak branch; the shear strain of the failure point was however maintained. The increase/decrease of the masonry Young modulus (MY) increased/decreased the initial stiffness and the first peak load, while a limited influence was detected in the post-peak.

**Table 12**. Main results of the sensitivity analysis.

| ID | Mod. parameter | Diagonal-compression | | | | In-plane bending | | | | Out-of-plane bending | | | |
|---|---|---|---|---|---|---|---|---|---|---|---|---|---|
| | | $F_{DCmax}$ [kN] | $\gamma(F_{DCmax})$ [‰] | $F_{DCu}$ [kN] | $\gamma(F_u)$ [%] | $F_{IBcr}$ [kN] | $\delta_{IBcr}$ [mm] | $F_{IBu}$ [kN] | $\delta_{IBu}$ [mm] | $F_{OBcr}$ [kN] | $\delta_{OBu}$ [mm] | $F_{OBu}$ [kN] | $\delta_{OBu}$ [mm] |
| Reference | | 350.5 | 0.25 | 308.9 | 0.93 | 36.2 | 0.06 | 80.5 | 5.5 | 33.1 | 0.4 | 41.7 | 34.1 |
| *GFRP yarns* | | | | | | | | | | | | | |
| YT +/- 33% | Peak/ultimate strains | 350.5 | 0.25 | 319.7 | 1.41 | 36.2 | 0.06 | 83.3 | 6.5 | 33.1 | 0.2 | 48.8 | 45.0 |
| | | 350.5 | 0.25 | 288.4 | 0.48 | 36.2 | 0.06 | 64.2 | 3.3 | 33.1 | 0.4 | 34.4 | 23.2 |
| YY +/- 33% | Young modulus | 350.5 | 0.25 | 324.2 | 0.74 | 36.2 | 0.06 | 82.5 | 4.1 | 33.1 | 0.4 | 43.9 | 26.2 |
| | | 350.5 | 0.25 | 285.4 | 1.29 | 36.2 | 0.06 | 77.4 | 7.7 | 33.1 | 0.4 | 42.5 | 44.6 |
| YA +/- 33% | Fiber area in yarns | 350.5 | 0.25 | 333.3 | 0.90 | 36.2 | 0.06 | 89.7 | 5.0 | 33.1 | 0.4 | 48.8 | 32.2 |
| | | 350.5 | 0.25 | 272.7 | 0.63 | 36.2 | 0.06 | 61.4 | 4.3 | 33.1 | 0.4 | 33.1 | 27.0 |
| YO | Orientation | 350.5 | 0.25 | 312.8 | 1.20 | 36.2 | 0.06 | 69.0 | 4.3 | 33.1 | 0.4 | 46.2 | 36.8 |
| *Mortar coating* | | | | | | | | | | | | | |
| CT +/- 33% | Tensile strength | 385.8 | 0.27 | 316.4 | 0.88 | 41.7 | 0.08 | 80.5 | 5.1 | 37.9 | 0.5 | 42.9 | 36.4 |
| | | 312.3 | 0.23 | 292.5 | 1.02 | 29.1 | 0.05 | 73.1 | 5.0 | 28.1 | 0.4 | 39.6 | 30.2 |
| CC +/- 33% | Compressive strength | 361.9 | 0.25 | 316.5 | 0.85 | 36.2 | 0.06 | 81.2 | 5.5 | 33.1 | 0.4 | 42.5 | 34.0 |
| | | 323.7 | 0.26 | 284.1 | 1.03 | 36.2 | 0.06 | 79.7 | 5.6 | 33.1 | 0.4 | 39.8 | 35.0 |
| CY +/- 33% | Young modulus | 335.5 | 0.22 | 304.8 | 0.93 | 38.0 | 0.06 | 81.9 | 5.6 | 35.1 | 0.4 | 42.9 | 34.5 |
| | | 371.2 | 0.30 | 294.4 | 0.89 | 32.7 | 0.06 | 78.1 | 6.0 | 31.0 | 0.4 | 39.7 | 32.4 |
| *Masonry* | | | | | | | | | | | | | |
| MT +/- 33% | Tensile strength | 385.4 | 0.33 | 310.3 | 0.89 | 41.5 | 0.07 | 80.8 | 5.6 | 35.4 | 0.6 | 42.2 | 35.8 |
| | | 295.3 | 0.26 | 305.5 | 0.98 | 28.5 | 0.05 | 80.7 | 5.5 | 28.9 | 0.4 | 41.2 | 33.6 |
| MC +/- 33% | Compressive strength | 352.5 | 0.25 | 339.8 | 0.92 | 36.2 | 0.06 | 80.5 | 5.5 | 33.1 | 0.4 | 41.8 | 33.8 |
| | | 342.2 | 0.26 | 262.1 | 0.97 | 36.2 | 0.06 | 74.9 | 5.8 | 33.1 | 0.4 | 38.5 | 36.4 |
| MY +/- 33% | Young modulus | 371.1 | 0.25 | 297.9 | 0.93 | 34.6 | 0.05 | 81.1 | 5.4 | 32.9 | 0.4 | 43.5 | 32.4 |
| | | 308.2 | 0.37 | 315.5 | 0.85 | 37.3 | 0.08 | 79.3 | 5.5 | 38.4 | 0.5 | 38.9 | 34.6 |

For in-plane and out-of-plane bending tests, samples IB-B-P66S-15 and OB-B-T66S-C8 were selected as reference; the main results are collected in Table 12. In both cases, the modification of the GFRP yarns parameters (Fig. 14b-c) influenced mainly the post-cracking behavior and led basically to the same conclusions of DC tests. The increasing of the tensile strength of mortar coating or of the masonry (CT, MT) led to higher first cracking loads and a more scattered post-peak branch; however, the coordinates of the failure point remained rather comparable. Modifications in the mortar Young modulus (CY) and in the compressive strength (CC) had a very slight influence.

Generally, for all the three test type, it was also observed that alterations (±33%) in the energy amount of the softening laws, in both tension (parameter $w_f/t$) and compression (parameter $a$), had a negligible influence for both the mortar coating and the masonry (<±5%).

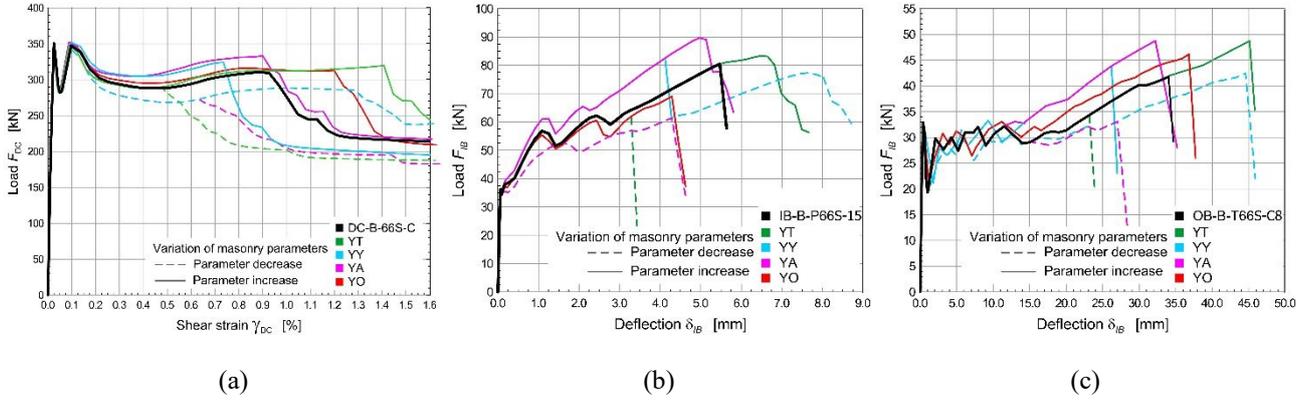

(a)            (b)            (c)

**Fig. 14** Sensitivity analysis: modification of GFRP yarn parameters in sample (a) DC-B-66S-C8, (b) IB-B-P66S-15 and (c) OB-B-T66S-C8.

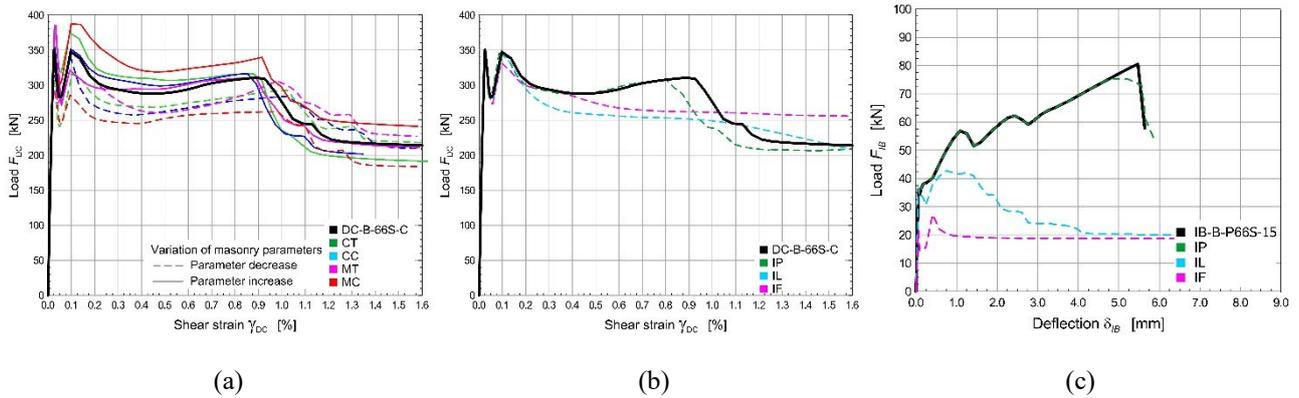

(a)            (b)            (c)

**Fig. 15** Sensitivity analysis: (a) influence of mortar/masonry resistance in DC-B-66S-C8 and role of the different interfaces in (b) DC-B-66S-C8 and (c) IB-B-P66S-15.

To investigate on the role of the connections, the assumption of no bonding was considered for the different interfaces, one by one (yarn-yarn, yarn-mortar and mortar-masonry).

When dealing with the DC-B-66S-C8 sample (Fig. 15b), it emerged that the absence of shear bond in the linear interfaces (IL), governing the yarn-mortar interaction, induced a visible resistance decrease in the post-peak branch of the $F_{DC}$-$\gamma_{DC}$ capacity curves. In fact, the slippage of the yarns in the most stressed area occurred as soon as the GFRP grid nodes resistance (point-to-point interfaces) was exceeded and did not allow the wires to exploit their whole resistance. Differently, the absence of shear bond in the point interfaces (IP), governing the yarn-yarn interaction, did not show any alteration in the curve trend, even though the failure point was anticipated. This difference is explained by considering the dominant role of the line interfaces against the yarn slippage, in respect to the less influent role of the point interfaces (as already emerged in the numerical study on pull-off tests performed in [7]). The no shear bond assumption at the mortar-masonry interface (IF) did not allow the collaboration between the two layers and, thus, the benefic damage diffusion of the reinforcement; therefore, the resistance decay occurred faster. It should be noted that, in actual application, this condition would also lead to the buckling of the mortar coating, not considered in the model.

Referring to the IB-B-P66S-15 sample, when the yarn-mortar bond was removed (IL) a rapid drop of resistance emerged as the yarn-yarn connections gradually failed along the lower yarn (Fig. 15c) and did not allow the yarn to attain to its ultimate strain. Differently, in the OB-B-T66S-C8 sample, the yarn failure was reached, even though the post peak branch exhibited a lower stiffness. For both IB and OB samples, no significant alterations were detected at the removal of the yarn-yarn connections (IR), just a slight anticipation of the failure point. The lack of any collaboration between mortar coating and the masonry (IF) led to the rapid activation of the rocking mechanism, as for the unreinforced masonry.

## 8. Conclusions

In this paper, a detail level modelling approach, using the OOFEM code, was applied to the simulation of CRM strengthened masonry elements under different loading conditions and also considering variations in the masonry and strengthening system characteristics. Recently, the model was independently calibrated on the basis of the experimental tests on materials and interfaces and validated through comparison with experimental characterization tests on CRM masonry samples (tensile, shear-bond and in-plane shear tests). The numerical study herein presented concerned diagonal compression tests, in-plane three-point bending tests and out-of-plane four-point bending tests of masonry elements. Experimental outcomes available in the literature were used for comparison, obtaining the validation of the model at a larger scale. The numerical model resulted capable to reproduce the typical failure modes of masonry elements (diagonal cracking, in-plane and out-of-plane bending) in terms of capacity curve and crack pattern; it was proved to be reliable also at advanced damage levels and at the varying of the characteristics of both the masonry and the reinforcement. The main discrepancies in the comparisons are related to intrinsic uncertainties in the material properties, possible cumulative damage due to the experimental loading-unloading procedure and a not precise "displacement control" in the experimental tests. A deep sensitivity analyses allowed to assess the main material parameters influencing the results. The findings of the analyses were somehow qualitatively expected but the reliability of the numerical model in catching such parametric variations was not obvious and needed to be assessed; also the quantitative influence of the parameters on the global performances was unknown and has been evaluated. Future research activities will concern the development and calibration of the simplified, intermediate level modelling, based on multilayered elements, to be adopted for the simulation of entire strengthened walls and structures. With this prospective, the detailed model, although very onerous and therefore difficult to apply on larger masonry elements, represents as useful tool to investigate on the validity of the assumptions of simplified models based, e.g., on perfect bond among the components.

**Data Availability Statement**

Some or all data, models, or code generated or used during the study are available in a repository online in accordance with funder data retention policies (https://doi.org/ 10.5281/zenodo.5827153). OOFEM v.2.5 was used for running the analyses (https://doi.org/10.5281/zenodo.4339630).

**Acknowledgments**

The project "conFiRMa" has received funding from the European Union's Horizon 2020 research and innovation programme under the Marie Sklodowska-Curie grant agreement No 101003410. The author wishes to thank the professional help of prof. B. Patzák, for in the use of the numerical tool, and the supervision of prof. A. Kohoutková (Czech Technical University in Prague - CTU).